%% file: dont-edit-this-main-m2.tex
\documentclass{amsart}
\usepackage[utf8]{inputenc}
\input{packages-and-definitions}

\title{Noetherian Operators in Macaulay2}
\author{Justin Chen, Yairon Cid-Ruiz, Marc H\"ark\"onen, Robert Krone, Anton Leykin}

\address{School of Mathematics, Georgia Institute of Technology,
Atlanta, Georgia}
\email{\{justin.chen,leykin\}@math.gatech.edu}
\email{harkonen@gatech.edu}

\address{Department of Mathematics: Algebra and Geometry, Ghent University, Krijgslaan 281 – S25, 9000 Gent, Belgium}
\email{Yairon.CidRuiz@UGent.be}

\address{Department of Mathematics, University of California, Davis, California}
\email{rkrone@math.ucdavis.edu}

\date{January 4, 2020}

\begin{document}
\input merge.tex

\maketitle

\begin{abstract}
A primary ideal in a polynomial ring can be described by the variety it defines and a finite set of Noetherian operators, which are differential operators with polynomial coefficients.
We implement both symbolic and numerical algorithms to produce such a description in various scenarios as well as routines for studying affine schemes through the prism of Noetherian operators and Macaulay dual spaces.
\end{abstract}

\section{Introduction}

The idea of describing ideals in polynomial rings via systems of differential operators was brewing since the beginning of the twentieth century.
In~\cite{Macaulay:modular-systems}, Macaulay brought forth the notion of an \emph{inverse system}, a system of differential conditions that describes a \emph{modular system} (a system of polynomials, or a polynomial ideal, in the modern language).

It was apparent to the contemporaries of Macaulay that a finite number of differential conditions should suffice to describe a $0$-dimensional affine or projective scheme.
In \cite{GROBNER_MATH_ANN}, Gr\"obner derived explicit characterizations for ideals that are primary to a rational maximal ideal~\cite[p.~174-178]{GROBNER_BOOK_AG_2}.
Moreover, he suggested that the same program can be carried out for any primary ideal~\cite[\S 1]{GROBNER_LIEGE}.

Despite this early algebraic interest, a complete description of primary ideals in terms of differential operators was first obtained by analysts in the \emph{Fundamental Principle of Ehrenpreis and Palamodov}~\cite{EHRENPREIS,PALAMODOV}.
At the core of the Fundamental Principle, one has the following theorem by Palamodov.

\begin{theorem}[{Palamodov}]
Let $R$ be a polynomial ring $R=\CC[x_1,\ldots,x_n]$ over the complex numbers, $P \subseteq R$ a prime ideal, and $Q \subseteq R$ a $P$-primary ideal.
There exist differential operators $A_1,\ldots,A_m \in R\langle \partial_{x_1},\ldots, \partial_{x_n}\rangle$ such that $Q = \lbrace f \in R \mid A_i \bullet f \in P \text{ for } 1 \le i \le m \rbrace$.
\end{theorem}

Following the terminology of Palamodov, the differential operators $A_1,\ldots,A_m$ are commonly called \emph{Noetherian operators for the $P$-primary ideal $Q$}.
Subsequent algebraic and computational approaches to characterize primary ideals with the use of differential operators have been given in \cite{BRUMFIEL_DIFF_PRIM}, \cite{OBERST_NOETH_OPS}, \cite{DAMIANO}, \cite{NOETH_OPS}; and, most recently, in \cite{chen2020noetherian} and \cite{cid2020primary}.

The purpose of this note is to present the \texttt{Macaulay2} package \texttt{NoetherianOperators}, which implements the algorithms for Noetherian operators introduced in \cite{chen2020noetherian} and \cite{cid2020primary} as well as the algorithms for (Macaulay) \emph{dual spaces} addressed in \cite{krone2017numerical,Krone-Leykin:eliminating-dual,Krone:dual-bases-for-pos-dim}.
While some of these algorithms rely on exact symbolic computation, the others employ numerical approximations using paradigms of \emph{numerical algebraic geometry}.

% \section{``Forward'' algorithms} 
\section{Computing Noetherian operators from ideals} \label{forward}

In this section we discuss four algorithms that can be used to compute a set of Noetherian operators for a primary ideal.
The main method in the \texttt{NoetherianOperators} package is \texttt{noetherianOperators}, which contains implementations of symbolic algorithms.
Specifying a value for the option \texttt{Strategy} allows the user to choose the algorithm that is used.
The method \texttt{numericalNoetherianOperators} implements a numerical algorithm, which can deal with approximate input.
Our list of strategies includes the following four algorithms:
\begin{enumerate}
    \item \textbf{Punctual Hilbert scheme}:
    the main idea of this algorithm is to use the \emph{punctual Hilbert scheme} to parametrize primary ideals.

    \item \label{list:symb} \textbf{Symbolic algorithm via dual spaces}: this algorithm computes Noetherian operators as bases of \emph{Macaulay dual spaces}.
    
    \item \label{list:nume} \textbf{Numerical algorithm via interpolation}: this algorithm interpolates Noetherian operators from their specializations at several general points sampled on the underlying variety.

    \item \textbf{Hybrid symbolic/numerical}: this approach optimizes the approach in (\ref{list:symb}) by using information obtained from applying the numerical algorithm (\ref{list:nume}) at one general point on the underlying variety.
\end{enumerate}

We now discuss each algorithm and illustrate its use in the package.
% divide the section into four short subsections that correspond to each one of the four algorithms (in the list above).
Throughout the article, let $\KK$ denote a field of characteristic zero, $R = \KK[x_1, \ldots, x_n]$ a polynomial ring over $\KK$, and $P \subseteq R$ a prime ideal in $R$.
We typically use $Q$ to denote a $P$-primary ideal, and $I$ to denote a general (not necessarily primary) ideal which has $P$ as a minimal prime.
% {\color{red} Talk about outputs, \texttt{DiffOp}?} 

To represent Noetherian operators, this package introduces a new type called \texttt{DiffOp}, representing elements in $R\langle \partial_{x_1}, \dotsc, \partial_{x_n} \rangle$.
A \texttt{DiffOp} is a hash table, where each key-value pair represents a term of the differential operator (with keys given by monomials in $R$, corresponding to a monomial in $\partial_{x_1}, \ldots, \partial_{x_n}$, whose value is the associated polynomial coefficient).
\texttt{DiffOp}s can be added and scaled, and also act on polynomials in $R$.
The output of \texttt{noetherianOperators} is a list of \texttt{DiffOp}s, and the output of \texttt{numericalNoetherianOperators} is a list of \texttt{InterpolatedDiffOp}s, a type similar to \texttt{DiffOp}, but whose coefficients are rational functions instead of polynomials.
\subsection{Punctual Hilbert scheme} \label{ssec:punctualHilbert}

The backbone of this algorithm is \cite[Theorem 2.1]{cid2020primary}, which can be seen as a ``representation theorem'' that parametrizes primary ideals via three 
% different but \justin{this will probably raise a comment from Dan Grayson}
closely related objects (points in the punctual Hilbert scheme, differentially closed vector spaces, and subbimodules of the Weyl-Noether module).

For a prime ideal $P$ of codimension $c$, let $\FF$ be the field of fractions of the integral domain $R/P$.
Up to a linear change of coordinates, we may (and do) assume that $\{ x_{c+1}, \ldots, x_n \}$ is a maximal independent set of variables modulo $P$, to simplify notation.

We now recall the steps of \cite[Algorithm 8.1]{cid2020primary}.
The main idea of this algorithm is to reduce the study of arbitrary $P$-primary ideals in $R$ to a zero-dimensional setting over the function field $\FF$.
This reduction is made by parametrizing $P$-primary ideals with the punctual Hilbert scheme
$
 {\rm Hilb}^m \bigl( \,\FF[[y_1,\ldots,y_c]] \,\bigr).
$
This is a quasiprojective scheme over the function field $\FF$.
Its classical points are ideals of colength $m$ in the local ring $\FF[[y_1,\ldots,y_c]]$.
For more details regarding punctual Hilbert schemes the reader is referred to \cite{IARROBINO}.
We define the inclusion map
\begin{equation*}
\gamma:R \hookrightarrow \FF[y_1,\dots,y_c]\, , \qquad
\begin{matrix}
x_i  &\mapsto &  y_i+\overline{x_i}, & \!\!\!\!\! \mbox{ for }1\leq i\leq c,\\
x_j & \mapsto  & \overline{x_j},& \quad \mbox{ for }c+1\leq j\leq n,
\end{matrix}
\end{equation*}
where $\overline{x_i}$ denotes the class of $x_i$ in $\FF$, for $1\leq i\leq n$.
With this, we can give the following explicit bijective correspondence
\begin{equation*}
\begin{array}{ccc}
\left\lbrace\begin{array}{c}
\mbox{$P$-primary ideals of $R$}\\
\mbox{with multiplicity $m$ over $P$}
\end{array}\right\rbrace
& \longleftrightarrow &
\left\lbrace\begin{array}{c}
\mbox{points in }{\rm Hilb}^m(\FF[[y_1,\ldots,y_c]])\\
\end{array}\right\rbrace\\
Q & \longrightarrow & J=\langle y_1,\dots,y_c\rangle^m+\gamma(
Q)\FF[y_1,\dots,y_c]\\
Q\,=\,\gamma^{-1}(J) & \longleftarrow & J.
\end{array}
\end{equation*}

The method \texttt{mapToPunctualHilbertScheme} can be used to compute the point in ${\rm Hilb}^m(\FF[[y_1,\ldots,y_c]])$ that corresponds to a $P$-primary ideal $Q$.
% Under the assumption of $\{ x_{c+1}, \ldots, x_n \}$ being a maximal independent set of variables modulo $P$,
For notational purposes, the method \texttt{mapToPunctualHilbertScheme} uses variables $\{hx_1,\ldots,hx_c\}$ (by adding an $h$ in front of each of $\{x_1,\ldots,x_c\}$) instead of  $\{y_1,\ldots,y_c\}$.
The following example illustrates the use of the method \texttt{mapToPunctualHilbertScheme}.

\beginOutput
i1 : needsPackage "NoetherianOperators";\\
\endOutput
\beginOutput
i2 : S = QQ[x_1, x_2, x_3];\\
\endOutput
\beginOutput
i3 : Q = ideal(x_1^2, x_2^2, x_1-x_2*x_3);\\
\endOutput
\beginOutput
i4 : mapToPunctualHilbertScheme Q\\
\                           2\\
o4 = ideal (hx  - x hx , hx )\\
\              1    3  2    2\\
\emptyLine
\                  /    S   {\char`\\}\\
o4 : Ideal of frac|--------|[hx ..hx ]\\
\                  |(x , x )|   1    2\\
\                  {\char`\\}  2   1 /\\
\endOutput
% \verbatiminput{generate_examples/punctual1.txt}
% \begin{verbatim}
%     i1 : needsPackage "NoetherianOperators";
%     i2 : R = QQ[x_1, x_2, x_3];
%     i3 : Q = ideal(x_1^2, x_2^2, x_1-x_2*x_3)
%                  2   2
%     o3 = ideal (x , x , - x x  + x )
%                  1   2     2 3    1
%     o3 : Ideal of R
%     i4 : mapIdealToPunctualHilbertScheme Q
%                                2
%     o4 = ideal (hx  - x hx , hx )
%                   1    3  2    2
%                       /    R   \
%     o4 : Ideal of frac|--------|[hx , hx ]
%                       |(x , x )|   1    2
%                       \  2   1 /
% \end{verbatim}

After computing the point $I \subseteq {\rm Hilb}^m(\FF[[y_1,\ldots,y_c]])$ corresponding to a $P$-primary ideal $Q$ of multiplicity $m$ over $P$, 
%the algorithm described in \cite[Algorithm 8.1]{cid2020primary} computes 
the inverse system $J^\perp$ of $J$ is computed.
% The last step in \cite[Algorithm 8.1]{cid2020primary} 
Lastly, an $\FF$-basis of $J^\perp$ is lifted to a set of Noetherian operators for the ideal $Q$.

% The method \texttt{noetherianOperators} when invoked with the option \texttt{Strategy => "PunctualHilbert"} uses the algorithm described 
% in \cite[Algorithm 8.1]{cid2020primary}
% above to compute a set of Noetherian operators.
This is the default strategy used to compute Noetherian operators, when the input is a primary ideal.
It can also be explicitly called by specifying \texttt{Strategy => "PunctualHilbert"}, as shown below:

\beginOutput
i5 : noetherianOperators(Q, Strategy => "PunctualHilbert")\\
o5 = \{1, x dx_1 + dx_2\}\\
\          3\\
\endOutput
% \verbatiminput{generate_examples/punctual2.txt}
% \begin{verbatim}
%     i5 : getNoetherianOperatorsHilb Q
%     o5 = {1, x dx  + dx }
%               3  1     2
%     o5 : set of Noetherian operators over the prime ideal (x , x )
%                                                             2   1
% \end{verbatim}

\subsection{Symbolic algorithm via dual spaces} \label{ssec:viaMacaulayMatrix}
The next algorithm to compute a set of Noetherian operators is a direct approach which reduces the problem to linear algebra.
% The main conceptual points in this approach are: (i) reduction to the zero-dimensional case, and (ii) a characterization of sets of Noetherian operators for zero-dimensional ideals, as bases of suitable dual spaces.
% To reduce to the zero-dimensional case, note that if $I$ is an $R$-ideal, and $\tt$ is a maximal set of independent variables on $R/I$, then the extension of $I$ to the ring $S := \KK(\tt)[\xx]$ is zero-dimensional, and $IS$ is primary if $I$ is primary.
% Furthermore, if $\{A_1, \ldots, A_m\} \subseteq R\langle \partial_{\tt}, \partial_{\xx} \rangle$ 
% {\color{red} Should define $W_A$ above} 
% are operators such that their images in $S\langle \partial_{\xx}\rangle$ are a set of Noetherian operators for $IS$, then $\{A_1, \ldots, A_m\}$ are in fact a set of Noetherian operators for $I$ (cf. \cite[Proposition 3.13]{chen2020noetherian}).
% Since one can always find ``lifts'' of any set of operators in $S\langle \partial_{\xx}\rangle$ (by clearing fractions), it suffices to compute a set of Noetherian operators for the zero-dimensional ideal $IS$, with the understanding that the underlying coefficient field $\KK(\tt)$ need not be algebraically closed.
For convenience, write $\tt$ for a maximal set of independent variables modulo $P$, and $\xx := \{ x_1, \ldots, x_n \} \setminus \tt$ as the dependent variables.
Then in the localization $S := \KK(\tt)[\xx]$ of $R$, the extension of $I$ to $S$ is zero-dimensional.
If now $I$ is a zero-dimensional $P$-primary ideal, then a set of Noetherian operators for $I$ is the same as a basis for the dual space of $I$ at $P$ (as will be discussed in in \Cref{sec:dualSpaces}).
This in turn can be computed as the kernel of a \emph{Macaulay matrix}, which is a matrix over $R/P$ with columns indexed by differential monomials and rows indexed by elements of $I$, whose entries are the result of applying a differential monomial to an element of $I$.
As the numbers of rows and columns increase, the kernel eventually stabilizes, at which point the result is returned.
% : namely, for a given degree $d$, let $f_1, \ldots, f_r$ be generators for $I$, and 
% \[
%   M_{d,I} := ( \partial^\beta \bullet (x^\alpha f_i) )_{|\alpha| < d, |\beta| \le d, 1 \le i \le r}
% \]
% be the matrix with entries in $R/P$ whose rows are indexed by polynomials $\{ x^\alpha f_i \mid |\alpha| < d, 1 \le i \le r \}$ and whose columns are indexed by differential monomials $\{ \partial^\beta \mid |\beta| \le d \}$.
% Then the kernel of $M_{d,I}$ is the truncation of the local dual space of $I$ at $P$ up to degree $d$, and for $d$ sufficiently large, a basis for this kernel will be a minimal set of Noetherian operators for $I$, cf. \cite[Algorithm 1]{chen2020noetherian}.

This is the default strategy used to compute Noetherian operators, when the input is a pair of ideals, the second of which should be a minimal prime of the first.
It can also be explicitly called by specifying \texttt{Strategy => "MacaulayMatrix"}, as shown below:
% The option \texttt{Strategy => "MacaulayMatrix"} will use this method to compute Noetherian operators.

\beginOutput
i6 : needsPackage "K3Carpets";\\
\endOutput
\beginOutput
i7 : I = carpet(2, 2, Characteristic => 0);\\
\emptyLine
o7 : Ideal of QQ[x ..x , y ..y ]\\
\                  0   2   0   2\\
\endOutput
\beginOutput
i8 : R = ring I;\\
\endOutput
\beginOutput
i9 : noetherianOperators(I, Strategy => "MacaulayMatrix")\\
o9 = \{1, 2y dy_0 + y dy_1\}\\
\           1        2\\
\endOutput

Calling \texttt{noetherianOperators} with a minimal prime ideal $P$ of $I$ as the second argument will compute a set of Noetherian operators for the $P$-primary component of $I$. If $I$ is unmixed, then the result of applying this method to all associated primes can be thought of as dual to a primary decomposition of $I$.

\beginOutput
i10 : (P1, P2) = (radical I, ideal(R_0 + R_1));\\
\endOutput
\beginOutput
i11 : J = intersect(I, P2^2);\\
\endOutput
\beginOutput
i12 : noetherianOperators(J, P1)\\
o12 = \{1, 2y dy_0 + y dy_1\}\\
\            1        2\\
\endOutput
\beginOutput
i13 : noetherianOperators(J, P2)\\
o13 = \{1, dx_0\}\\
\endOutput

\subsection{Numerical algorithm via interpolation} \label{ssec:numericalAlg}
We also provide algorithms to compute Noetherian operators purely from numerical data, bypassing the need to compute Gr\"obner bases.
This is based on computing a set of \emph{specialized Noetherian operators}, i.e. the result of evaluating (at some point) all polynomial coefficients in a set of Noetherian operators.
The key observation is that one can obtain a set of specialized Noetherian operators by suitably slicing the variety \cite[Theorem 4.1]{chen2020noetherian}.
More precisely, for a $P$-primary ideal $Q \subseteq \CC[\tt,\xx]$ and a point $p = (\tt_0, \xx_0) \in V(P)$, a minimal set of specialized Noetherian operators corresponds to a basis of the dual space of the zero-dimensional ideal $Q + (\tt - \tt_0)$ at the point $p$.
The function \texttt{specializedNoetherianOperators} can be used to perform this computation.

%\verbatiminput{generate_examples/numerical1.txt}
\beginOutput
i14 : p = point\{\{1,1,1,1,1,1\}\};\\
\endOutput
\beginOutput
i15 : specializedNoetherianOperators(I, p, DependentSet => \{R_1, R_3, R_4\})\\
o15 = \{1, 2*dy_0 + dy_1\}\\
\endOutput

Once a set of specialized Noetherian operators has been computed at a single general point, subsequent computations at other points can be sped up as the monomial support of a valid set of Noetherian operators is known (this fact also underlies the hybrid approach in \Cref{ssec:hybrid}).
After specialized Noetherian operators are computed at sufficiently many general points on the variety, the original set of Noetherian operators can be recovered from their specializations via interpolation of rational functions, cf. \cite[Algorithm 5]{chen2020noetherian}.

This is the preferred strategy when the input is inexact, although it can also be used for exact input, as shown below. 
Note that the value of \texttt{DependentSet} must be specified:

%\verbatiminput{generate_examples/numerical2.txt}
\beginOutput
i16 : numericalNoetherianOperators(I, DependentSet => \{R_1, R_3, R_4\})\\
\                    x\\
\                     1\\
o16 = \{1, 1*dy_0 + ---*dy_1\}\\
\                   2x\\
\                     0\\
\endOutput

By default, \texttt{Bertini} is used to sample points on $V(\sqrt{I})$. The user can specify their own sampling function with the option \texttt{Sampler}. The sampler should be a function that takes an integer $n$ and the ideal $I$ as input, and return a \texttt{List} of $n$ points on the variety.
\beginOutput
i17 : needsPackage"NumericalImplicitization";\\
\endOutput
\beginOutput
i18 : numericalNoetherianOperators(I, DependentSet => \{R_1, R_3, R_4\},\\
Sampler => (n,I) -> apply(n, i -> point sub(matrix realPoint radical I, CC)))\\
\                    x\\
\                     1\\
o18 = \{1, 1*dy_0 + ---*dy_1\}\\
\                   2x\\
\                     0\\
\endOutput

\subsection{Hybrid symbolic/numerical} \label{ssec:hybrid}
In \Cref{ssec:viaMacaulayMatrix}, Noetherian operators are found by computing the kernel of a Macaulay matrix with entries in the function field $\KK(\tt)$, but computations in this field can be expensive.
The numerical approach in \Cref{ssec:numericalAlg} instead specializes the independent variables to random values.
This allows computations to be done in $\KK$ (typically with $\KK = \CC$) which is cheaper but the result consists of specializations of the Noetherian operators.  %An interpolation step is needed to recover the operators themselves.
A hybrid approach can combine the best of both strategies.
% In this algorithm, specializations of Noetherian operators $\{A_1(p),\ldots,A_m(p)\}$ for $I$ are found at a random point $p$ as in \Cref{ssec:numericalAlg}.
% Computing $A_i(p)$ over $\KK$ reveals the support $C$ of $A_i$ in the monomials $\partial_{\xx}$.
% Then $A_i$ is computed over $\KK(\tt)$ from only the columns of the Macaulay matrix indexed by $C$.
% This reduces the size of the computation needed in the rational function field. 
In essence, the information revealed from running the numerical algorithm at a single point (without performing interpolation) can be used to trim the Macaulay matrix down to a smaller (optimal) size, without changing the kernel.
For more details, we refer the interested reader to \cite[Section 4.1]{chen2020noetherian}.

This strategy is called by specifying \texttt{Strategy => "Hybrid"} with the method \texttt{noetherianOperators}. On larger examples, this strategy can greatly outperform the approach in \Cref{ssec:viaMacaulayMatrix}.

\beginOutput
i19 : noetherianOperators(I, Strategy=>"Hybrid")\\
o19 = \{1, 2x y dy_0 + x y dy_1\}\\
\            0 2        2 1\\
\endOutput

As in the numerical algorithm, the user may also specify a sampling function to find a general point.
% This allows the computation of Noetherian operators for components of unmixed ideals.

% \section{``Backward'' algorithm}
\section{Computing ideals from Noetherian operators}

In this section, we discuss a procedure that can be seen as the inverse of the process of computing a set of Noetherian operators. 
% Let $\KK$ be a field of characteristic zero, $R=\KK[x_1,\ldots,x_n]$ be a polynomial ring over $\KK$, and $P \in \Spec(R)$ be a prime ideal. 
First, note that for any $R$-bisubmodule $\mathcal{E} \subseteq R\langle \partial_{x_1},\ldots, \partial_{x_n}\rangle$ of the Weyl algebra, the set 
$$
\lbrace f \in R \mid A \bullet f \in P \text{ for all } A \in \mathcal{E} \rbrace
$$
is always a $P$-primary ideal in $R$ (cf. \cite[Proposition 3.5]{NOETH_OPS}). We now consider the following problem: 

\begin{center} \label{prob:backwardsProblem}
Given an $R$-bisubmodule $\mathcal{E} \subseteq R\langle \partial_{x_1},\ldots, \partial_{x_n}\rangle$, compute (generators for) the $P$-primary ideal 
$\lbrace f \in R \mid A \bullet f \in P \text{ for all } A \in \mathcal{E} \rbrace$.
\end{center}

This is accomplished by \cite[Algorithm 8.2]{cid2020primary}.
The idea is to use the explicit maps provided in \cite[Theorem 2.1]{cid2020primary} in inverse order to how they appear in \cite[Algorithm 8.1]{cid2020primary} (i.e., as discussed in \Cref{ssec:punctualHilbert}).
It should be noted that our implementation solves the following 
% slightly more general 
effective version of the problem above:

\begin{center}
Given $A_1,\ldots,A_m \in R\langle \partial_{x_1},\ldots, \partial_{x_n}\rangle$, compute the $P$-primary ideal 
$\lbrace f \in R \mid A \bullet f \in P \text{ for all } A \in \mathcal{E} \rbrace$, \\
where $\mathcal{E} \subseteq R\langle \partial_{x_1},\ldots, \partial_{x_n}\rangle$ is the $R$-bisubmodule generated by $A_1,\ldots,A_m$. 
\end{center}

The function \texttt{getIdealFromNoetherianOperators} implements \cite[Algorithm 8.2]{cid2020primary}. Below we show an example in which given a $P$-primary ideal $Q$, we compute a set of Noetherian operators for $Q$, and then we recover $Q$ from its Noetherian operators along with $P$. In general, this process may result in a different set of generators for $Q$.

\beginOutput
i20 : R = QQ[x_1,x_2,x_3];\\
\endOutput
\beginOutput
i21 : Q = ideal(x_1^2, x_2^2, x_3^2, x_1*x_2 + x_1*x_3 +x_2*x_3);\\
\endOutput
\beginOutput
i22 : L = noetherianOperators Q\\
o22 = \{1, dx_3, dx_2, dx_1, - dx_1*dx_2 + dx_2*dx_3, - dx_1*dx_2 + dx_1*dx_3\}\\
\endOutput
\beginOutput
i23 : Q' = getIdealFromNoetherianOperators(L, radical Q)\\
\              2   2                       2\\
o23 = ideal (x , x , x x  + x x  + x x , x )\\
\              3   2   1 2    1 3    2 3   1\\
\endOutput
\beginOutput
i24 : Q == Q'\\
o24 = true\\
\endOutput

\section{Dual spaces and local Hilbert functions} \label{sec:dualSpaces}

In Section \ref{forward}, dual spaces were used in service of computing Noetherian operators.
However, dual spaces can also directly provide information about polynomial ideals.
Suppose $P$ is the maximal ideal corresponding to a $\KK$-rational point $p \in \KK^n$, and $I \subseteq R$ an ideal with $p \in V(I)$.
The \emph{dual space of $I$ at $P$} is
\[ D_P[I] := \{A \in \KK[\partial_{x_1}, \ldots, \partial_{x_n}] \mid (A \bullet f)(p) = 0 \text{ for all } f \in I\}. \]
The dual space is a subspace of the space of differential operators on $R$ with constant coefficients and finite support which uniquely determines $IR_P$, where $R_P$ denotes the localization of $R$ at $P$.

The following dual space algorithms work with numerical data, for example if $\KK = \CC$ and the point associated to $P$ is known only approximately.
This is in contrast to symbolic algorithms relying on Gr\"obner bases.
The methods described in this section were previously part of a package titled \texttt{NumericalHilbert} which has now been incorporated into \texttt{NoetherianOperators} due to an overlap in functionality.

If $P$ is a minimal prime of $I$, then the dual space is finite dimensional, and a basis of the dual space is a set of Noetherian operators for the $P$-primary component of $I$.
Otherwise the dual space is infinite dimensional, and only a truncation up to a specified degree can be computed.
The methods \texttt{zeroDimensionalDual} and \texttt{truncatedDual} compute a basis for the dual space in these two cases.
As in \Cref{ssec:viaMacaulayMatrix}, these dual spaces are computed as kernels of Macaulay matrices.
From a basis of the dual space truncated to degree $d$, it is straightforward to compute the local Hilbert function of $R/I$ up to degree $d$, and this is implemented as applying \texttt{hilbertFunction} to a \texttt{DualSpace} object.

\beginOutput
i25 : R = CC[x_1,x_2];\\
\endOutput
\beginOutput
i26 : I = ideal\{x_1^2 + x_2^2 - 4, (x_1 - 1)^2\};\\
\endOutput
\beginOutput
i27 : p = point\{\{1.0, 1.7320508\}\};\\
\endOutput
\beginOutput
i28 : D = zeroDimensionalDual(p, I)\\
o28 = | 1 -1.73205x_1+x_2 |\\
o28 : DualSpace\\
\endOutput
\beginOutput
i29 : apply(3, i -> hilbertFunction(i, D))\\
o29 = \{1, 1, 0\}\\
o29 : List\\
\endOutput

%{\color{red} What about a Hilbert function computation here?}

Another way of truncating the dual space of a positive dimensional ideal is with respect to a local elimination order instead of a degree order.
In \cite{Krone-Leykin:eliminating-dual} this is referred to as an \emph{eliminating dual space}.
Assume $\{x_{c+1},\ldots,x_n\}$ is a maximal set of independent variables for $R/I$, and choose an order that eliminates $V = \{x_1,\ldots,x_c\}$.
An eliminating dual for $I$ with respect to $V$ truncated to degree $d$, denoted $E^d_P[I,V]$, is defined as the set of dual operators whose lead terms with respect to the monomial order have degree at most $d$ in the variables $V$, and is computed with method \texttt{eliminatingDual}.
When $V = \{x_1\}$, such a truncated dual space allows one to find the dual space of the colon ideal $I:x_1$ directly (without requiring the potentially expensive symbolic computation of finding $I:x_1$):
 \[ E_P^{d}[I:x_1, \{x_1\}] = x_1 \cdot E_0^{d+1}[I,\{x_1\}] \]
 where $x_1 \cdot A$ represents the right action of $x_1 \in R$ on $A \in \KK[\partial_{x_1}, \ldots, \partial_{x_n}]$ (for example $x_1 \cdot \partial_{x_1}^2 = 2\partial_{x_1}$).
A representation of these colon ideals are needed in \cite[Algorithm 5.1]{Krone-Leykin:eliminating-dual} for identifying embedded primes on a curve.
%Eliminating dual spaces also have a connection to Noetherian operators for positive dimensional ideals. {\color{red} Should this be explained?}

Using \cite[Algorithm 23]{Krone:dual-bases-for-pos-dim}, computing truncated dual spaces up to a certain degree provides a numerical algorithm for finding a full set of generators for the initial ideal of $I$ with respect to a local degree order, and this is implemented by \texttt{gCorners}.
In the process, a standard basis can be computed by specifying \texttt{ProduceSB => true}.
\beginOutput
i30 : gCorners(point p, I)\\
o30 = | x_2 x_1^2 |\\
\              1       2\\
o30 : Matrix R  <--- R\\
\endOutput
Finding the initial ideal via approximate numerical methods is an essential part of \texttt{isPointEmbedded}, which implements a numerical algorithm for the detection of an embedded component developed in \cite[Algorithm 4.2]{krone2017numerical}.

\bibliography{references}{}
\bibliographystyle{plain}

\end{document}

%% file: packages-and-definitions.tex
\usepackage[utf8]{inputenc}
\usepackage{amsfonts,amssymb,amsbsy,amsmath} % Math
\usepackage{bm} % Bold
\usepackage[inline,shortlabels]{enumitem} % Inline enumeration
\usepackage{listings} % Code snippets
\usepackage{amsthm} % Theorems
\usepackage{commath} % Common math
\usepackage{booktabs} % Tables
\usepackage[letterpaper, margin=1in]{geometry}
\usepackage{relsize}
\usepackage{tikz-cd}
\usepackage{tikz}
\usepackage{booktabs}
\usepackage{verbatim}
\usepackage{algorithm,algpseudocode}
\usepackage{todonotes}
\usepackage{stmaryrd}
\usepackage{mathtools}
\usepackage{hyperref}
\usepackage{cleveref} % Clever reference
\usepackage{kbordermatrix}
\usetikzlibrary{patterns}

\theoremstyle{plain}
\newtheorem{theorem}{Theorem}[section]

\theoremstyle{definition}

\newcommand{\B}[1]{\mathbb #1}

\newcommand{\CC}{{\B C}}

\newcommand{\KK}{{\B K}}

\newcommand{\FF}{{\B F}}

\newcommand{\xx}{\boldsymbol{x}}
\renewcommand{\tt}{\boldsymbol{t}}

\algrenewcommand\algorithmicrequire{\textbf{Input}}
\algrenewcommand\algorithmicensure{\textbf{Output}}

\allowdisplaybreaks

%% file: merge.tex
% these macros are needed by the output of the 'merge' program:

{\obeyspaces\global\let =\ }
\font\tteight=cmtt8
\font\ttten=cmtt10

\newdimen\outputBaseLineSkip
         \outputBaseLineSkip = 10pt % was 9pt
\newskip\beginOutputSkip
        \beginOutputSkip = 4.5 pt plus .5 pt minus .5 pt
\newskip\endOutputSkip
        \endOutputSkip   = 2.5 pt plus .25 pt minus .25 pt

\def\looserOutput#1{%
  \advance\beginOutputSkip by #1
  \advance\endOutputSkip by #1
}
\def\tighterOutput#1{%
  \advance\beginOutputSkip by -#1
  \advance\endOutputSkip by -#1
}

\def\beginOutput{%
    \par
    \penalty -150
    % \vskip \beginOutputSkip
    \penalty -150
    \begingroup
      \def\\{%
          \leavevmode
          \hss
          \endgraf
          \penalty 150
          }
      % \tteight
      \ttten
      % \baselineskip = \outputBaseLineSkip
      \parindent = 24pt
      \def\${\char`\$}
      \def\{{\char`\{}
      \def\}{\char`\}}
      \catcode`\_=\the\catcode`a
      \catcode`\^=\the\catcode`a
      \catcode`\#=\the\catcode`a
      \catcode`\~=\the\catcode`a
      \catcode`\&=\the\catcode`a
      \parskip=0pt
      \lineskip=0pt
      \obeyspaces
      }
\def\emptyLine{%
    \penalty -100
    % \vskip \endOutputSkip
    % \vskip \beginOutputSkip
    \penalty -100
    }
\def\endOutput{%
    \endgroup
    \par
    \penalty -150
    % \vskip \endOutputSkip
    \penalty -150
    \noindent}